\input amstex
\documentstyle{amsppt}
\magnification=\magstep1
\NoRunningHeads
%\NoBlackBoxes
\topmatter
\title
 New spectral multiplicities~for~mixing transformations
\endtitle
\author  Alexandre~I.~Danilenko
\endauthor
%\rightheadtext{}
\address
Institute for Low Temperature Physics \& Engineering
of Ukrainian National Academy of Sciences,
47 Lenin Ave.,
 Kharkov, 61164, UKRAINE
\endaddress
\email            alexandre.danilenko\@gmail.com
\endemail

\thanks
\endthanks
%\subjclass Primary 28D15, 28D20; Secondary 46L55
%\endsubjclass
\abstract
It is shown that if $E$ is any subset of $\Bbb N$ such that  either $1\in E$ or $2\in E$  then there is a {\it mixing} transformation whose set of spectral multiplicities is $E$.
\endabstract

 \endtopmatter
 \document
%\NoBlackBoxes

\head 0. Introduction
\endhead

Let $T$ be an invertible  measure preserving transformation of a standard probability non-atomic space $(X,\goth B,\mu)$.
Denote by $U_T$ the associated unitary Koopman operator on the space $L^2(X,\mu)$:
$$
U_Tf=f\circ T.
$$
Let $\Cal M(T)$ stand for the set of essential values of the spectral multiplicity function for $U_T\restriction (L^2(X,\mu)\ominus\Bbb C)$.
A subset $M\subset\Bbb N$ is called {\it realizable} if $M=\Cal M(T)$ for an ergodic transformation $T$.
The {\it spectral multiplicity problem} can be stated as follows
\roster
\item"{\bf (Pr1)}" which subsets of $\Bbb N$ are realizable?
\endroster
It is one of the basic problems of the spectral theory of dynamical systems.
Many papers are devoted to partial solutions of (Pr1): \cite{Os}, \cite{Ro1}, \cite{Ro3}, \cite{G--Li}, \cite{KL}, \cite{Ag1}, \cite{Ry1}, \cite{Ag2}, \cite{Ag3}, \cite{Da3}, \cite{Ag4}, \cite{KaL}, \cite{Da5},
\cite{Ry5} (see also a recent survey \cite{Le}).
In particular, it was shown there that the following  subsets are realizable:
\roster
\item"---" every subset containing $1$,
\item"---"
every subset containing $2$,
\item"---"
subsets $n\cdot E$, where $n>0$ and $E$ is an arbitrary subset containing $1$.
\endroster
Despite of this progress, (Pr1) remains unsolved in general.
A transformation $T$ is said to be
\roster
\item"$\circ$"{\it mixing} if $\mu(T^nA\cap B)\to\mu(A)\mu(B)$ as $n\to\infty$ for all subsets $A,B\in\goth B$,
\item"$\circ$"{\it rigid} if $\mu(T^{n_k}A\cap B)\to\mu(A\cap B)$ for some subsequence $n_k\to\infty$ and  all subsets $A,B\in\goth B$,
\item"$\circ$"{\it partially rigid} if $\mu(T^{n_k}A\cap B)\ge\delta\mu(A\cap B)$ along some subsequence $n_k\to\infty$ for certain $\delta>0$ and  all subsets $A,B\in\goth B$.
\endroster
We note that the realizations in the aforementioned papers are either rigid or partially rigid. The present paper is devoted to a more subtle question:
 \roster
\item"{\bf (Pr2)}" for which subsets $E$ of $\Bbb N$ there is a {\it mixing} transformation $T$ such that $\Cal M(T)=E$?
\endroster
We then say that $E$ admits a {\it mixing realization}.
The information about mixing realizations is rather scarce.
The first mixing transformation with a simple spectrum (i.e. a mixing realization for $\{1\}$) was constructed
in \cite{Ne} as a Gaussian system.
In \cite{Ro2} it was  shown that for each $n>0$, there is a mixing $T$ with $n<\sup\Cal M(T)<\infty$.
The paper \cite{Ry3} provides mixing realizations for $\{2\}$.
More generally, for each $n>1$, Ageev builds mixing realizations for $\{2,3\dots,n\}$ in \cite{Ag4}.
It  should be noted that his result follows also from \cite{Ry4}, where Ryzhikov constructs mixing transformations whose  symmetric powers  have all a simple spectrum.
(We also mention that (P1) and (P2) are  solved completely in the category of {\it infinite} measure preserving transformations \cite{DaR1}, \cite{DaR2}).

The set of mixing transformations $\Bbb M$ is meager in the (Polish) weak topology on the group of all measure preserving transformations.
Tikhonov \cite{Ti} introduced a Polish topology on $\Bbb M$ which is strictly stronger then the induced weak topology. He showed, in particular, that a generic mixing transformation  has a simple spectrum.
Thus we are interested in  realizations belonging to a meager part of $\Bbb M$.
This somehow explains a complexity of (Pr2) comparatively with (Pr1): meager sets in Tikhonov's topology are ``much smaller'' then meager sets in the weak topology.

We now state the main result of the paper.

\proclaim{Main Theorem}
Given an arbitrary  subset  $E$ of \, $\Bbb N$
such that either $1\in E$ or $2\in E$,
there exists  a mixing (of all orders) transformation $S$ such that $\Cal M(S)=E$.
\endproclaim

To prove the theorem we combine various techniques and approaches to spectral realizations developed by a number of authors:
\roster
\item"$\bullet$"
non-Abelian   compact group extensions \cite{Ro1}--\cite{Ro3},
\item"$\bullet$"
natural factors of Abelian compact group extensions \cite{KL}, \cite{G--L},
\item"$\bullet$"
spectral multiplicities of Cartesian products and their natural factors
 \cite{Ka}, \cite{Ry1}, \cite{Ag1}, \cite{Ag4},
\item"$\bullet$"
weak limits of powers of the associated Koopman operator \cite{Ag1},
\cite{Ry1}, \cite{KaL}, etc.,
\item"$\bullet$"
 techniques to force  mixing  \cite{Ry4},
\item"$\bullet$"
mixing properties of staircase algorithms \cite{Ad},
\item"$\bullet$"
$(C,F)$-construction \cite{dJ}, \cite{Da2}--\cite{Da5}.
 \endroster

The outline of the paper is as follows. In Section~1 we consider cocyles of dynamical systems with values in compact group extensions and discuss their simplest spectral properties.
$(C,F)$-construction of dynamical systems, an algebraic version of the classical cutting-and-stacking  is  outlined briefly in Section~2.

In Section~3 we prove  a part  of Main Theorem when $1\in E$.
It is worthy to note that the   realizations from \cite{G--Li} and \cite{KL} based on Abelian compact group extensions  are---on our opinion---inherently rigid and hence do not suit for this purpose.
That is why we ``return'' to  historically earlier non-Abelian group extensions that go back to \cite{Ro1}, \cite{Ro3} (see also \cite{KaL}).
 Since they were elaborated  for spectral realization of only {\it finite} sets, we first show how to adjust them  to {\it arbitrary} subsets
of $\Bbb N$ via a ``compactification'' of the Algebraic Lemma from \cite{KL}.
Then we adapt Ryzhikov's idea   {\it to force mixing} in symmetric products (\cite{Ry3} and \cite{Ry4}) to  compact extensions. The idea is as follows.
A $(C,F)$-construction of the extensions is implemented in such a way to ensure availability of some special weak limits in the closure of powers of the associated Koopman operators.
These limits imply simple spectrum of some invariant components  of the Koopman operators.
That, in turn, guarantees that the spectral multiplicities can be computed from a certain algebraic picture.
Moreover, for each $\delta>0$, it is possible to modify this construction in such a way to obtain a ``mixing component''  that ``occupies'' $\delta$-part of the whole space. At the same time the ``rigid component'' occupying $(1-\delta)$-part of the space is retained to guarantee the necessary weak limits of powers of the Koopman operator.
Then we {\it pass to a limit} along a sequence of such $\delta_i$-``partly mixing'' transformations which have the same (desired) spectral multiplicities as $\delta_i\to 1$. The limit will be a mixing compact extension of a rank-one mixing map.
Though there are no non-trivial weak limits in the closure of powers of the limit Koopman operator any more, it is possible to retain the {\it simple spectrum} properties of the components of the Koopman operators in the limit. Hence the limit Koopman operator has the same spectral multiplicities as the prelimit ones.

 Section~4 is devoted entirely to the proof of the remaining part  of  Main Theorem, i.e. when $2\in E$. The scheme of the proof is as in the previous section.
 We note that the realizations (of all sets containing $2$) constructed in \cite{Da5} are immanently  rigid. They are  not suit to force mixing.
Therefore we again use non-Abelian group extensions. For that we first
generalize the Algebraic Lemma from \cite{KL} in such a way to {\it
embrace} the case when $1\notin E$. It should be noted that our
generalization implies the affirmative answer to the following question
from \cite{KaL}:
 \roster
\item"" {\it whether the series of weak mixing realizations constructed
    in  \cite{KaL} covers all finite subsets of $\Bbb N$ that contain $2$?}
 \endroster
Several concluding remarks are contained in Section~5.

The author thanks V.~Ryzhikov    for  a suggestion  to apply the techniques   ``to force mixing'' from \cite{Ry3} and \cite{Ry4}  to compact group extensions.

\head 1. Cocycles with vales in group extensions
\endhead
Let $T$ be an ergodic transformation of $(X,\mu)$. Denote by $\Cal R\subset X\times X$ the $T$-orbit equivalence relation. A Borel map $\alpha$ from $\Cal R$ to a compact group $K$ is called a {\it cocycle} of $\Cal R$ if
$$
\alpha(x,y)\alpha(y,z)=\alpha(x,z)
 $$
for all $(x,y), (y,z)\in\Cal R$. Two cocycles $\alpha,\beta:\Cal R\to K$ are {\it cohomologous} if
$$
\alpha(x,y)=\phi(x)\beta(x,y)\phi(y)^{-1}
$$
at a.a. $(x,y)\in \Cal R$ for a Borel map $\phi:X\to K$. Denote by $C(T)$ the {\it centralizer} of $T$, i.e. the group of all transformations commuting with $T$. Given  $S\in C(T)$, a cocycle $\alpha\circ S:\Cal R\to K$ is well defined by
$$
\alpha\circ S(x,y):=\alpha(Sx,Sy).
 $$
for all $(x,y)\in \Cal R$.
Let $\lambda_K$ denote Haar measure on $K$. Define a transformation
$T_\alpha$ of the probability space $(X\times K,\mu\times\lambda_G)$ by setting
$$
T_\alpha(x,k):=(Tx,k\alpha(x,Tx)).
$$
It is called the $\alpha$-{\it skew product extension of} $T$. If it is ergodic then $\alpha$ is called {\it ergodic}.
We note that the centralizer of $T_\alpha$ contains a compact subgroup $K$ which acts on $X\times K$ by left translations along the second coordinate.
We note that if $K$ is Abelian then $U_{T_\alpha}$ is unitarily equivalent to the orthogonal sum
$
\bigoplus_{\chi\in\widehat K}U_{T,\chi},
$
where $U_{T,\chi}$ is a unitary operator in $L^2(X,\mu)$ given by
$$
U_{T,\chi}f(x)=\chi(\alpha(x,Tx))f(Tx).
$$
If $H$ is a closed subgroup of $K$ then we denote by $T_{\alpha,H}$  the following transformation
of $(X\times H\backslash K,\mu\times\lambda_{H\backslash K})$:
$$
T_{\alpha,H}(x,Hk):=(Tx,Hk\alpha(x,Tx)).
$$
It is easy to see that $T_{\alpha,H}$ is a factor of $T_{\alpha}$. It is called a {\it natural} factor of $T_{\alpha}$.
Notice that $U_{T_\alpha,H}$ is unitarily equivalent to $\bigoplus_{\chi\in\widehat{K/H}}U_{T,\chi}$.

Let   $A$ be a compact $K$-module. This means that $A$ is a compact second countable Abelian group  and $K$ acts on $A$ by group automorphisms.
Denote by $G$ the semidirect product $K\ltimes A$. Recall that the multiplication in $G$ is given by
$$
(k,a)(k',a'):=(kk', a+k\cdot a').
$$
Given a cocycle $\gamma:\Cal R\to G$, we  can write it as a pair of ``coordinate'' maps
$\gamma=(\beta,\alpha)$. Then $\beta$ is a cocycle of $\Cal R$ with values in $K$. As for  the map $\alpha:\Cal R\to A$, it is not a cocycle but it satisfies the following equation
$$
\alpha(x,y)+\beta(x,y)\cdot \alpha(y,z)=\alpha(x,z)
$$
for all $(x,y),(y,z)\in\Cal R$.
Denote by $\Cal R(\beta)$ the $T_\beta$-orbit equivalence relation on $X\times K$. It is easy to verify that
$$
\Cal R(\beta)=\{((x,k), (x',k\beta(x,x')))\mid (x,x')\in\Cal R,k\in K\}.
$$
Consider a map
$$
\widetilde\alpha:\Cal R(\beta)\ni ((x,k), (x',k\beta(x,x')))
\mapsto k\cdot\alpha(x,x')\in A.
$$
 Then $\widetilde\alpha$ is a cocycle of $\Cal R(\beta)$. Moreover,
 the $\gamma$-skew product extension $T_\gamma$ coincides with the ``double'' extension  $(T_\beta)_{\widetilde\alpha}$ of $T$.
Given $k\in K$, denote by $L_k$ the following transformation of $(X\times K,\mu\times\lambda_K)$: $L_k(x,r):=(x,kr)$. Then $L_k\in C(T_\beta)$ and
$$
\widetilde\alpha\circ L_k=k\cdot\widetilde\alpha.
$$
This  implies the following:
$$
U_{T_\beta,\chi}\quad\text{is unitarily equivalent to $U_{T_\beta,\chi\circ k}$ for each }k\in K.
\tag1-1
$$
For more information about cocycles  with values in locally compact group extensions we refer to \cite{Da1} and references therein.

\head 2. $(C,F)$-systems and $(C,F)$-cocycles
\endhead

We recall the  $(C,F)$-construction (see \cite{dJ}, \cite{Da2}--\cite{Da4}).
Let two sequences $(C_n)_{n>0}$ and $(F_n)_{n\ge 0}$ of finite subsets in $\Bbb Z$ are given such that:
\roster
\item"---"
$F_n=\{0,1,\dots,h_n-1\}$,
\item"---"  $0\in C_n$, $\# C_n>1$,
\item"---" $F_n+C_{n+1}\subset F_{n+1}$,
\item"---" $(F_{n}+c)\cap (F_n+c')=\emptyset$ if $c\ne c'$, $c,c'\in C_{n+1}$,
\item"---" $\lim_{n\to\infty}h_n/(\#C_1\cdots\# C_n)<\infty$.
\endroster
Let $X_n:=F_n\times C_{n+1}\times C_{n+2}\times\cdots$. Endow this set with the (compact Polish) product topology. The  following map
$$
(f_n,c_{n+1},c_{n+2},\dots)\mapsto(f_n+c_{n+1},c_{n+2},\dots)
$$
is a topological embedding of $X_n$ into $X_{n+1}$. We now set $X:=\bigcup_{n\ge 0} X_n$ and endow it with the (locally compact Polish) inductive limit topology. Given $A\subset F_n$, we denote by $[A]_n$ the following cylinder: $\{x=(f,c_{n+1},\dots,)\in X_n\mid f\in A\}$. Then $\{[A]_n\mid A\subset F_n, n>0\}$ is the family of all compact open subsets in $X$. It forms a base of the topology on $X$.

Let  $\Cal R$ stand for the {\it tail} equivalence relation on $X$: two points $x,x'\in X$ are $\Cal R$-equivalent if there is $n>0$ such that $x=(f_n,c_{n+1},\dots),\ x'=(f_n',c_{n+1}',\dots)\in X_n$ and $c_m=c_m'$ for all $m>n$. There is only one probability (non-atomic) Borel measure $\mu$ on $X$ which is invariant (and ergodic) under $\Cal R$.

Now we define a transformation $T$ of $(X,\mu)$ by setting
$$
T(f_n,c_{n+1},\dots):=(1+f_n,c_{n+1},\dots )\text{ whenever }f_n<h_n-1,\ n>0.
$$
This formula defines $T$ partly on $X_n$. When $n\to\infty$, $T$ extends to the entire $X$ (modulo a countable subset) as a $\mu$-preserving invertible transformation. Moreover, the $T$-orbit equivalence relation  coincides  with $\Cal R$ on the subset where $T$ is defined. We call $(X,\mu,T)$ {\it the $(C,F)$-dynamical system associated with $(C_{n+1},F_n)_{n\ge 0}$}.

From now on we will assume that $\# C_{n}\to\infty$ as $n\to\infty$.

 If for each $n>0$, there are $i_n\ge 0$ and  $c_n\in C_n$  such that $c_n+ih_n+\frac{i(i-1)}{2}\in C_n$ for each $i=1,\dots,i_n$ and $i_n/\# C_n\to 1$ as $n\to\infty$ then the corresponding $(C,F)$-system is called an {\it almost staircase} \cite{Ry4}.
It follows from \cite{Ad}  that if $\frac{i_n^2}{h_n}\to 0$ then the corresponding almost staircase transformation is mixing.

We recall a concept of $(C,F)$-cocycle (see \cite{Da3}). Let $G$ be a compact second countable group.  Given a sequence of maps $\gamma_n:C_n\to G$, $n=1,2,\dots$, we first define a Borel cocycle $\gamma:\Cal R\cap (X_0\times X_0)\to G$ by setting
$$
\gamma(x,x'):=\gamma_1(c_1)\cdots\gamma_n(c_n)\gamma_n(c_n')^{-1}\cdots \gamma_1(c_1')^{-1},\tag2-1
$$
if $x=(f_0,c_1,c_2,\dots)\in X_0$,  $x'=(f_0',c_1',c_2',\dots)\in X_0$ and $c_i=c_i'$ for all $i>n$. To extend $\gamma$ to the entire $\Cal R$, we first define a map $\pi:X\to X_0$ as follows. Given $x\in X$, let $n$ be the least positive integer such that $x\in X_n$. Then $x=(f_n,c_{n+1},\dots)\in X_n$. We  set
$$
\pi(x):=(\underbrace{0,\dots,0}_{n+1\text{ times}}, c_{n+1}, c_{n+2},\dots)\in X_0.
$$
Of course, $(x,\pi(x))\in\Cal R$. Now for each pair $(x,y)\in\Cal R$, we let
$$
\gamma(x,y):=\gamma(\pi(x),\pi(y)).
$$
It is easy to verify that $\gamma$ is a well defined cocycle of $\Cal R$ with values in $K$. We call it {\it the $(C,F)$-cocycle associated with} $(\gamma_n)_{n=1}^\infty$.

Suppose now that $G=K\ltimes A$. Then we can write $\gamma=(\beta,\alpha)$ as above. Also,  $\gamma_n=(\beta_n,\alpha_n)$ for some maps $\beta_n:C_n\to K$ and $\alpha_n:C_n\to A$. Then $\beta$ is the $(C,F)$-cocycle associated with $(\beta_n)_{n>0}$. As for the mapping $\alpha$, we deduce from~\thetag{2-1} that
$$
\aligned
\alpha(x,x')&=\alpha_1(c_1)+\sum_{j=2}^{n}\beta_1(c_1)\cdots
\beta_{j-1}(c_{j-1})\cdot\alpha_j(c_j)\\
&-\beta(x,x')\cdot\bigg[\alpha_1(c_1')+\sum_{j=2}^{n}\beta_1(c_1')\cdots
\beta_{j-1}(c_{j-1}')\cdot\alpha_j(c_j')\bigg]
\endaligned
\tag2-2
$$
for all $x=(f,c_1,c_2,\dots),\ x'=(f',c_1',c_2',\dots)\in X_0$ such that
$c_i=c_i'$ whenever $i>n$.

\head 3. Mixing transformations~whose~set~of~spectral multiplicities contains 1
\endhead
Let $E$ be an arbitrary  subset of \, $\Bbb N$
such that $1\in E$.
Our main purpose in this section is to prove the following theorem.

\proclaim{Theorem 3.1} There exists  a mixing (of all orders) transformation $S$ such that $\Cal M(S)=E$.
\endproclaim

Let $K$ be a compact Abelian group, $B$ a countable $K$-module and $D$ a subgroup of $B$. We denote
$$
L(K,B,D):=\{\#(\{k\cdot d\mid k\in K\}\cap D), d\in D\setminus\{0\}\}.
$$

The  lemma below plays an important role in the proof of Theorem~3.1.
 We show how to deduce it  from  the Algebraic Lemma from \cite{KL}.

\proclaim{Lemma 3.2}
There  are a compact  Abelian group $K$, a compact $K$-module $A$
and a  closed subgroup $H\subset A$ such that $L(K,\widehat A,\widehat{A/H})=E$ and
\roster
\item"\rom{(i)}"
 $K$ is the inverse limit of a sequence of finite  Abelian groups $K_1\leftarrow K_2\leftarrow\cdots$, where the arrows are onto homomorphisms,
\item"\rom{(ii)}"
 $A$  is the inverse limit of  a sequence of finite $K$-modules $A_1\leftarrow A_2\leftarrow\cdots$, where the arrows are $K$-equivariant onto homomorphisms (more precisely, $A_n$ is a $K_n$-module, i.e.
$k\cdot a=k'\cdot a$ if $\pi_n(k)=\pi_n(k')$, where $\pi_n:K\to K_n$ stands for the canonical projection),
\item"\rom{(iii)}"
there is a dense countable sub-$K$-module $A_0\subset A$ such that the $K$-orbit of each $a\in A_0$ is finite.
\endroster
\endproclaim

\demo{Proof}
It is shown in \cite{KL} that there is a compact group $A=\Bbb T^{k_1}\times\Bbb T^{k_2}\times\cdots$,  a closed subgroup $H$ of $A$ and
a group automorphism $\sigma=\sigma_1\times\sigma_2\cdots$ of $A$ such that $L(\Bbb Z, \widehat A, \widehat{A/H})=E$, where $\sigma_i$ is the automorphism of $\Bbb T^{k_i}$ generated by the cyclic permutation of the coordinates and $\widehat A$ is considered as a $\Bbb Z$-module with respect to the automorphism which is dual to $\sigma$.
Our first observation is that the argument from
\cite{KL} works as well if we replace
every copy of
$\Bbb T$  in $A$ with $\Bbb Z/2\Bbb Z$.
Next, it is easy to see that the closure of the cyclic group generated by $\sigma$ in the group of all group automorphisms of $A$ equipped  with the natural Polish topology is the inverse limit of the sequence
$$
\Bbb Z/k_1\Bbb Z\leftarrow \Bbb Z/(k_1k_2\Bbb Z)\leftarrow\cdots
$$
equipped with the natural compact topology.
Denote this inverse limit by $K$.
 Then $L(\Bbb Z, \widehat A, \widehat{A/H})=L(K,\widehat A,\widehat{A/H})$.
The properties (i)--(iii) are now obvious.
\qed
\enddemo

From now on till the end of the section $(K,A,H)$ denote the groups whose existence is asserted in Lemma~3.2.

We first construct   a partially rigid  transformation $S$ with
 $\Cal M(S)=E$ and whose  {\it staircase (i.e. mixing) part} is $(1-\delta)$-large.
For that we will use skew product extensions of $(C,F)$-transforma\-tions constructed via $(C,F)$-cocycles.

Fix a countable dense subgroup $K_0\subset K$.
Let $0<i_n\le r_n$, $r_n\to\infty$,  $i_n/r_n\to\delta$, as $n\to\infty$, where $0<\delta< 1$.
We now let
$$
c_{n+1}(i):=
\cases
0, &\text{ if } i=0,\\
c_{n+1}(i-1)+h_n, &\text{ if } 0<i< i_n,\\
c_{n+1}(i-1)+h_n+i-i_n, &\text{ if } i_n\le i< r_n.
\endcases
$$
We define $C_{n+1}$
by setting
$$
C_{n+1}:=\{c_{n+1}(i)\mid 0\le i<r_n\}.
$$
Let   $A_0\subset A$ be a countable submodule which  is dense in $A$.
 Consider a partition $\Bbb N=\bigsqcup_{a\in A_0}\Cal N_a\sqcup
\bigsqcup_{k\in K_0}\Cal L_k$ of $\Bbb N$ into infinite subsets $\Cal N_a$, $\Cal L_k$.

If $n+1\in\Cal N_a$, we define a mapping $\alpha_{n+1}: C_{n+1}\to A$ by
$$
\alpha_{n+1}(c_{n+1}(i)):=
\cases
0, &\text{ if } i=0,\\
\alpha_{n+1}(c_{n+1}(i-1))- a, &\text{ if }i=1,\dots,i_n-1,\\
\alpha_{n+1}(c_{n+1}(i-1)), &\text{ if }i=i_n,\dots,r_n-1.
\endcases
$$
We also let $\beta_{n+1}\equiv 1$ on $C_{n+1}$.

If $n+1\in\Cal L_k$, we let $\alpha_{n+1}\equiv 0$ on $C_{n+1}$.
 We also define $\beta_{n+1}: C_{n+1}\to K$  by setting
$$
\beta_{n+1}(c_{n+1}(i)):=
\cases
0, &\text{ if } i=0,\\
\beta_{n+1}(c_{n+1}(i-1))- k, &\text{ if }i=1,\dots,i_n-1,\\
\beta_{n+1}(c_{n+1}(i-1)), &\text{ if }i=i_n,\dots,r_n-1.
\endcases
$$

Let $(X,\mu,T)$ stand for  the $(C,F)$-system associated with $(C_{n+1},F_n)_{n\ge 0}$. Denote by $\Cal R$ the tail equivalence relation on $X$ and denote by $\gamma=(\beta,\alpha):\Cal R\to G$ the $(C,F)$-cocycle associated with $(\beta_n,\alpha_n)_{n>0}$.
Recall that given a character $\chi\in\widehat A$,  $U_{T_\gamma,\chi}$ stand for the following unitary operator on $L^2(X\times K,\mu\times\lambda_K)$:
$$
U_{T_\beta,\chi}v(x,k):=\chi(\widetilde\alpha((x,k),(Tx,k\beta(x,Tx))))
v(Tx,k\beta(x,Tx)).
$$
We also note that  $U_{T_\beta, 1}:=U_{T_\beta}$.
Given $a\in A_0$, we  let $\text{Orb}_K(a):=\{k\cdot a\mid k\in K\}$ and
$$
l_\chi(a):=\frac 1{\# \text{Orb}_K(a)}\sum_{b\in \text{Orb}_K(a)}\chi(b)
$$
for $\chi\in\widehat A$ and $a\in A$.

\proclaim{Proposition 3.4} For each $a\in A_0$,
\roster
\item"\rom{(i)}" $U^{h_n}_{T_\beta,1}\to\delta \,I+(1-\delta)P_0$ weakly as $\Cal N_a-1\ni n\to\infty$, where $P_0$ is the orthogonal projection onto the subspace of constant functions on $X\times K$.
\item"\rom{(ii)}"
$U^{h_n}_{T_\beta,\chi}\to\delta \, l_\chi(a)I$ weakly as $\Cal N_a-1\ni n\to\infty$ for each $1\ne \chi\in\widehat A$.
\endroster
For each $k\in K_0$,
\roster
\item"\rom{(iii)}" $U^{h_n}_{T,1}\to\delta \,I+(1-\delta)P_0$ weakly as $\Cal L_k-1\ni n\to\infty$, where $P_0$ is the orthogonal projection onto the subspace of constant functions on $X$.
\item"\rom{(iv)}"
$U^{h_n}_{T,\eta}\to\delta \, \eta(k)I$ weakly as $\Cal L_k-1\ni n\to\infty$ for each $1\ne \eta\in\widehat K$.
\endroster
\endproclaim

\demo{Proof}
We first verify that $\beta$ is ergodic.
Take arbitrary $n>0$, $k\in K$  and a subset $B\subset F_n$.
Find $m>n$ such that $m\in\Cal L_k-1$.
Let
$$
I_{m}:= \{0 \le i<i_m-1\mid \}.
$$
We now let
$$
B_0:=B+C_{n+1}+\cdots+C_m+c_{m+1}(I_m)\subset F_{m+1}.
$$
If $x=(f_m,c_{m+1}(i_1),c_{m+2}(i_2),\dots)\in X_m$ belongs to $[B_0]_{m+1}$
then $i_1\in I_m$ and
$$
T^{h_m}x=(f_m,c_{m+1}(i_1+1), c_{m+2}(i_2),\dots).
$$
Moreover,
$$
[B_0]_{m+1}\cup T^{h_m}[B_0]_{m+1}\subset [B]_n,\  \mu([B_0]_{m+1})>0.5\delta\mu([B]_n),
$$
and $\beta(x,T^{h_m}x)=k$ for all $x\in [B_0]_{m+1}$.
By a standard criterium of ergodicity for cocycles (see \cite{FM}, \cite{Sc}),  $\beta$ is ergodic.

Our next purpose is to show (i).
 Take $n\in\Cal N_a-1$,  subsets $B,B'\subset F_n$ and $\eta\in\widehat K$. Then
$$
\aligned
(U^{h_n}_{T,\eta} 1_{[B]_n})(x)
&=
\eta(\beta(x,T^{h_n}x))1_{T^{-h_n}[B]_n}(x)\\
&=\sum_{j=1}^{i_n-1}\eta(\beta(x,T^{h_n}x))1_{[B+c_{n+1}(j-1)]_{n+1}}(x)\\
&+
\sum_{j=i_n+1}^{r_n-1} \eta(\beta(T^{i_n-j}x,T^{h_n}x))  U^{i_n-j}_{T_\beta,\eta} 1_{[B+c_{n+1}(j-1)]_{n+1}}(x) +\bar o(1)\\
&=\sum_{j=1}^{i_n-1}1_{[B+c_{n+1}(j-1)]_{n+1}}(x)\\
&+
\sum_{j=i_n+1}^{r_n-1}  U^{i_n-j}_{T,\eta} 1_{[B+c_{n+1}(j-1)]_{n+1}}(x) +\bar o(1)
\endaligned
$$
for all $x\in X$ and $k\in K$.
Then
$$
\langle U^{h_n}_{T,\eta}1_{[B]_n},1_{[B']_n}\rangle =
\frac{i_n}{r_n}\langle 1_{[B]_n}, 1_{[B']_n}\rangle
+\bigg\langle\frac 1{r_n}\sum_{j=1}^{r_n-i_n}
 U^{-j}_{T,\eta} 1_{[B]_n},1_{[B']_n}\bigg\rangle+\bar o(1).
$$
It follows from the mean ergodic theorem that
$(r_n-i_n)^{-1}\sum_{j=1}^{r_n-i_n}
 U^{-j}_{T,\eta}$ converges strongly to the orthogonal projection $P_\eta$ onto the subspace of $U_{T,\eta}$-invariant vectors.
Hence
$$
U^{h_n}_{T,\eta}\to \delta I+(1-\delta)P_\eta \text{ \ weakly as }
\Cal N_a-1\ni n\to\infty.
$$
If $\eta=1$ then $P_\eta=P_0$ since $T$ is ergodic.
If $\eta\ne 1$ and $U_{T_\beta,\eta} v=v$ for a non-trivial function $v\in L^2(X,\mu)$ then $\eta(\beta(x,Tx))v(Tx)=v(x)$ at a.a. $x$.
It follows that $\eta\circ\beta$ is cohomologous to the trivial cocycle.
However this contradicts to the fact that $\beta$ is ergodic.
Hence $v=0$ and therefore $P_\eta=0$.
Since $U_{T_\beta,1}=U_{T_\beta}=\bigoplus_{\eta\in\widehat K}U_{T,\eta}$,
the claim~(i) follows.

Claims (iii) and (iv) are shown  in a similar way.

 It remains to prove (ii). Fix $n\in \Cal N_a-1$.
Take $B\subset F_n$, $ \eta\in\widehat K$ and $\chi\in\widehat A$.
It follows that
$$
\aligned
&U^{h_n}_{T_\beta,\chi}(1_{[B]_n}\otimes\eta)(x,k)
=
\chi(k\cdot\alpha(x,T^{h_n}x))1_{T^{-h_n}[B]_n}(x) \eta(k\beta(x,T^{h_n}x))\\
&=\sum_{j=1}^{i_n-1}\chi(k\cdot\alpha(x,T^{h_n}x))
1_{[B+c_{n+1}(j-1)]_{n+1}}(x) \eta(k)\\
&+
\sum_{j=i_n+1}^{r_n-1} \eta(\beta(T^{i_n-j}x,T^{h_n}x))  U^{i_n-j}_{T_\beta,\chi}(1_{[B+c_{n+1}(j-1)]_{n+1}}\otimes\eta)(x,k) +\bar o(1)
\endaligned
\tag3-1
$$
for all $x\in X$ and $k\in K$.

We say that two points $x,y\in X$ are $\Cal R_n$-equivalent if either $x=y$ or $x,y\in X_n$ and the ``tails'' of $x,y$ in $X_n=F_n\times C_{n+1}\times\cdots$ coincide starting from the second coordinate.
If $x=(f_n,c_{n+1}(i_{n+1}),\dots)\in X_n$ then we let $\bar x^n:=(f_n,0,0,\dots)\in X_n$. It now follows from \thetag{2-2} that
$$
\alpha(x,T^{h_n}x)=\beta(\bar x^n,0)\cdot a
\quad\text{for all  }x\in X_n':=\bigsqcup_{0\le i<i_n-1}[F_n+c_{n+1}(i)]_{n+1}.
$$
For each $b\in \text{Orb}_K(a)$, we let
$$
Y_{n,b}:=\bigg\{(x,k)\in X_n'\times K\,\bigg|\, k\cdot\alpha(x,T^{h_n}x)=b\bigg\}.
$$
It is easy to see that
\roster
\item"$(\heartsuit)$"
$\{Y_{n,b}\mid b\in \text{Orb}_K(a)\}$
is a partition of  $X_n'\times K$ into subsets of equal measure.
\endroster
Moreover,
$$
\multline
\sum_{j=1}^{i_n-1}\chi(k\cdot\alpha(x,T^{h_n}x))
1_{[B+c_{n+1}(j-1)]_{n+1}}(x)\eta(k)\\
=
\sum_{b\in \text{Orb}_K(a)}\chi(b)1_{Y_{n,b}}(x,k)
1_{[B]_{n}}(x)\eta(k)+\bar o(1).
\endmultline
\tag3-2
$$
We claim that for each $b\in \text{Orb}_K(a)$,
\roster
\item"$(\lozenge)$"
the set $Y_{n,b}$ is
$\Cal R_n(\beta)$-invariant,
\endroster
i.e. if $(x,k)\in Y_{n,b}$ and
$(x,y)\in\Cal R_n$ then $(y,k\beta(x,y))\in Y_{n,b}$.
This follows from the following fact:
$$
\beta(\bar x^n,0)=\beta(\bar x^n,\bar y^n)\beta(\bar y^n,0)=
\beta(x,y)\beta(\bar y^n,0)
$$
whenever $(x,y)\in \Cal R_n$.

 Since $\bigcup_{n>0}\Cal R_n(\beta)=\Cal R(\beta)$  and $\beta$ is ergodic and  the unit ball of bounded linear operators in $L^2(X\times K,\mu\times\lambda_K)$ is weakly compact  and the algebra
$L^\infty(X\times K)$ is weakly closed and $\mu(X_n')\to\delta$, we deduce from
$(\lozenge)$ and $(\heartsuit)$ that $1_{Y_{n,b}}$ weakly converges  to
$\delta(\#\text{Orb}_K(a))^{-1}I$
for each $b\in\text{Orb}_K(a)$  as $\Cal N_a-1\ni n\to\infty$.
Hence
$$
\sum_{b\in \text{Orb}_K(a)}\chi(b)1_{Y_{n,b}}\to \delta l_\chi(a)I \text{ \ weakly as }\Cal N_a-1\ni n\to\infty.
\tag3-3
$$

We now show that the cocycle $\widetilde\alpha$ (or, equivalently, $\gamma$) is ergodic.
Take arbitrary $n>0$,  $a\in A_0$,  subset $B\subset F_n$ and
an open subset $U\subset K$.
Select $m>n$ such that $m\in\Cal N_{a}-1$ and
$$
\mu\times\lambda_K(([B]_n\times U)\cap Y_{m,a})>\frac \delta{2\# \text{Orb}_K(a)}
(\mu\times\lambda_K)([B]_n\times U).
$$
If $(x,k)\in ([B]_n\times U)\cap Y_{m,a}$
then $T^{h_m}_\beta(x,k)\in [B]_n\times U$ and
$$
\widetilde\alpha((x,k), T_\beta^{h_m}(x,k))=k\cdot\alpha(x,T^{h_m}x)=a.
$$
 By the  standard criterium of ergodicity for cocycles,  the $A_j$-projection of $\widetilde\alpha$ is ergodic for every $j>0$.
It follows that $\widetilde\alpha$ (or $\gamma$) is ergodic.

Finally,  we deduce (ii) from  \thetag{3-1} by applying  \thetag{3-3} and \thetag{3-2}.
To see  that the second sum in  \thetag{3-1}  vanishes in  the limit one should repeat  the corresponding argument from  the proof of (i) almost literally and make use of the ergodicity of~$\gamma$ (instead of $\beta$).
\qed
\enddemo

While proving the above proposition we established as a byproduct that $\gamma$ is ergodic. We now show a stronger fact.

\proclaim{Corollary 3.5}
The $\gamma$-skew product extension $T_\gamma$ of $T$ is weakly mixing.
\endproclaim

\demo{Proof} If $T_\gamma v=\lambda v$ for a vector $0\ne v\in L^2_0(X\times G,\mu\times\lambda_G)$ and some $\lambda\in\Bbb T$ then
$|\langle U_{T_\gamma}^{q}v,v\rangle|=\|v\|^2$ for each $q>0$. On the other hand, it follows from Proposition~3.4(i) and (ii) that $|\langle U_{T_\gamma}^{h_n}v,v\rangle|\to \delta\|v\|^2$ as $\Cal N_0-1\ni n\to\infty$, a contradiction.
\qed
\enddemo

Given $b\in A$ and $k\in K$, let $L_b$ denote the operator of translation by $b$ along the $A$-coordinate in
$L^2(X\times K\times A,\mu\times\lambda_K\times\lambda_A)$ and
let $L_k$  denote the operator of translation by $k$ along the $K$-coordinate in
$L^2(X\times K,\mu\times\lambda_K)$.

\remark{Remark \rom{3.6}} The assertion of Proposition 3.4 can also be restated in the following equivalent way:
$$
\align
U_{T_\gamma}^{h_n}&\to
\delta\cdot\frac{\sum_{b\in
\text{Orb}_K(a)}L_{b}}{\#\text{Orb}_K(a)}+(1-\delta)P_0 \text{ \ weakly
as }\Cal N_a-1\ni n\to\infty \quad\text{and}\\
U_{T_\beta}^{h_n}&\to \delta L_k+(1-\delta)P_0 \text{ \ weakly
as }\Cal L_k-1\ni n\to\infty\\
\endalign
$$
for each $a\in A_0$ and $k\in K_0$.
\endremark

Given $\chi,\chi'\in\widehat A$, we write $\chi\sim\chi'$ if there is $k\in K$ with $\chi=\chi'\circ k$.

\proclaim{Theorem 3.7}
 $\Cal M((T_{\beta})_{\widetilde\alpha,H})=L(K,\widehat A,\widehat{A/H})=E$.
\endproclaim

\demo{Proof}
Since $U_{T_{\gamma,H}}=\bigoplus_{\chi\in \widehat{A/H}}U_{T_\beta,\chi}$,
it is enough to show the following:
\roster
\item"(a)"
$U_{T_\beta,\chi}$ has a simple spectrum for each $\chi\in\widehat A$,
\item"(b)"
$U_{T_\beta,\chi}$ and $U_{T_\beta,\chi'}$
are unitarily equivalent if $\chi'\sim\chi$,
\item"(c)"
the measures of maximal spectral type of $U_{T_\beta,\chi}$ and $U_{T_\beta,\chi'}$ are mutually singular if $\chi'\not\sim\chi$.
\endroster

We first note that (b) follows from \thetag{1-1}.

Now, if $\chi'\ne\chi\circ k$ for any $k\in K$ then there is $a\in A_0$ such that
$l_\chi(a)\ne l_{\chi'}(a)$.
This follows from (ii) and (iii) of Lemma~3.2.
Pass to the   weak limits in the sequences $U^{h_n}_{T_\gamma,\chi}$ and
$U^{h_n}_{T_\beta,\chi'}$ as $\Cal N_a-1\ni n\to\infty$.
Then (c) follows from  (i) and (ii) of Proposition~3.4.

It remains to show (a).
For this purpose we  decompose $U_{T,\beta}$ into a direct orthogonal sum
$$
U_{T_\beta}=\bigoplus_{\eta\in\widehat K}U_{T,\eta}.
$$
Since $T$ is of rank one and the function $X\ni x\mapsto \eta(\beta(x,Tx))\in\Bbb T$ is constant on each cylinder $[f]_{n}$, $0\le f\le h_{n}-2$, it follows that the operator
$U_{T,\eta}$ has a simple spectrum.
We also deduce from (iii) and (iv) of Proposition~3.4 that the maximal spectral types of $U_{T,\eta}$ and $U_{T,\eta'}$ are mutually disjoint if $\eta\ne\eta'$.
Hence $U_{T_\beta}$ has a simple spectrum.
Moreover, it is easy to verify that
\roster
\item"---"
for each $n$,  vectors $1_{[f]_n}\otimes 1_{\pi^{-1}_n(k)}$, $0\le f\le h_{n}-2$, $k\in K_n$, belong to the $U_{T_\beta}$-cyclic space generated by $\sum_{\eta\in\widehat K_n} 1_{[0]_n}\otimes\eta$  (at least up to small $\epsilon$ in the $L^2$-norm) and
\item"---"
the function
$$
X\times K\ni (x,k)\mapsto(\chi(k\cdot\alpha(x,Tx)),\beta(x,Tx))\in\Bbb T^2
$$
is constant on each `level' $[f]_n\times \pi^{-1}_n(k)$, $0\le f<h_n-1$, $k\in K_n$.
\endroster
These two facts imply (a).
\qed
\enddemo

We note that the  condition (a) and (c) together  imply
$$
U_{T_\beta,\chi}\oplus U_{T_\beta,\xi} \text{ has a simple spectrum for all $\xi,\chi\in\widehat A$, $\chi\not\sim\xi$}.\tag 3-5
$$
In turn, \thetag{3-5} implies
\roster
\item"---"
$T_{\gamma}$  is weakly mixing and
\item"---"
$\Cal M(T_{\gamma,H})=E$.
\endroster

To force mixing we fix a sequence of positive reals $\delta_j$ that goes to $0$ very fast and construct for each $j>0$,
\roster
\item"$\bullet$" a $(C,F)$-system
$(X^{(j)},\mu_j, T_j)$ associated to some sequence
$(C_n^{(j)},F_{n-1}^{(j)})_{n>0}$ and
\item"$\bullet$" a cocycle
$\gamma^{(j)}=(\beta^{(j)},\alpha^{(j)}):\Cal R_j\to G=K\ltimes A$
associated to  some sequence $(\beta_n^{(j)},\alpha_n^{(j)})_{n>0}$
\endroster
such that
 \thetag{3-5} is satisfied for the pair $(T_j,\gamma^{(j)})$ and
the set $C_{n+1}^{(j)}$ is the union of an arithmetic sequence
$\{0,h_n^{(j)},\dots,(i_n^{(j)}-1)h_n^{(j)}\}$ with a {\it staircase} \, sequence
$$
\bigg\{i_n^{(j)}h_n^{(j)}+1,\dots, (r_n^{(j)}-1)h_n^{(j)}
+\frac{(r_n^{(j)}-i_n^{(j)})(r_n^{(j)}- i_n^{(j)}+1)}{2}\bigg\}
$$
and $i_n^{(j)}/r_n^{(j)}\to\delta_j$ as $n\to\infty$.
Of course,  $\Cal R_j$  stands here for the $T_j$-orbit equivalence relation on $X^{(j)}$.
In what follows we call such $(T_j,\gamma^{(j)})$  a $\delta_j$-{\it pair}.
The pair to be found will appear as a certain limit of the sequence $(T_j,\gamma^{(j)})$. Our main concern is to preserve \thetag{3-5} in the limit.

In the beginning of this section we outlined an explicit algorithm how to construct $\delta$-pairs for an arbitrary $\delta>0$, an arbitrary {\it initial} set $F_0$ and an arbitrary sequence $r_n\to\infty$ (provided that $\mu(X)<\infty$).

We carry out the construction of the sequence
$(T_j,\gamma^{(j)})$ via an inductive procedure in $j$.

Let $(T_1,\gamma^{(1)})$ be an arbitrary $\delta_1$-pair.
Suppose now that we have already built $(T_j,\gamma^{(j)})$.
It follows from the inverse limit structure of $K$  that
$$
l^2(K_1)\subset l^2(K_2)\subset\cdots
 $$
 and $\bigcup_{m>0}l^2(K_m)$  is dense in $L^2(K,\lambda_K)$.
Denote by $\Cal H_s^{(j)}$    the finite dimensional linear space $\bigoplus_{f\in F_s^{(j)}}1_{[f]_s}\otimes l^2(K_s)$.
Then
$\Cal H_1^{(j)}\subset\Cal H_2^{(j)}\subset\cdots$ and $\bigcup_{s>0}\Cal H_s^{(j)}$
is dense in $L^2(X_j\times K,\mu_j\times\lambda_K)$.
Since \thetag{3-5} is satisfied for $(T_j,\gamma^{(j)})$,
 there exist  $s_j>0$, $Q_j>0$ and  a vector    $v_j\in\Cal H_{s_j}^{(j)}\oplus\Cal H_{s_j}^{(j)}$ such that
\roster
\item"$\bullet$"
for each $w\in \Cal H_{s_j}^{(j)}\oplus\Cal H_{s_j}^{(j)}$,
$$
\inf_{a_{-q},\dots,a_q\in\Bbb R}\bigg\| w-\sum_{|q|\le Q_j}a_{q}
\bigg(U^q_{(T_j)_{\beta^{(j)}},\chi}\oplus U^q_{(T_j)_{\beta^{(j)}},\xi}\bigg)v_j\bigg\|_j\le \delta_j
\|w\|_j
$$
whenever  $\chi\not\sim\xi\in\widehat A_j$ and
\item"$\bullet$"
$U^q_{(T_j)_{\gamma^{(j)}},\chi}v\in\Cal H_{s_j}^{(j)}$ for all $|q|\le Q$ and $\chi\in\widehat A_j$,
\endroster
where $\|.\|_j$ denotes the norm in $L^2(X^{(j)}\times K,\mu_j\times\lambda_K)$.
Now let $(T_{j+1},\gamma_{j+1})$ be any $\delta_{j+1}$-pair whose initial set  $F_0^{(j+1)}$ is $F_{s_j}^{(j)}$ and  the {\it sequence of cuts} $(r_n^{(j+1)})_{n>0}$ is $j+1,j+2,\dots$.
Thus, we now think that the sequence $(C_{n+1}^{(j+1)},F_n^{(j+1)})_{n\ge 0}$ is constructed.
It remains to repeat  this process infinitely many times.

We now  construct a certain ``$(C,F)$-limit'' for the sequence
$(T_j,\gamma^{(j)})_{j=1}^\infty$.
For that we just let
 $$
\align
s_0&:=0,\\
F_n &:=
F^{(j)}_{n- s_1-\cdots-s_{j-1}} \quad\text{ if }      s_0+\cdots+s_{j-1}\le n<s_1+\cdots+s_j,\\
C_n &:=C^{(j)}_{n- s_1-\cdots-s_{j-1}} \quad\text{ if }      s_0+\cdots+s_{j-1}< n\le s_1+\cdots+s_j \quad\text{ and }\\
\gamma_n &:=\gamma^{(j)}_{n- s_1-\cdots-s_{j-1}} \quad\text{ if }       s_0+\cdots+s_{j-1}< n\le s_1+\cdots+s_j.
\endalign
$$
Denote by $(X,\mu,T)$ the $(C,F)$-system associated to $(C_n,F_{n-1})_{n>0}$.
Let $\gamma$ stand for the cocycle associated to $(\gamma_n)_{n>0}$.
We now show that $T_{\gamma,H}$ {\it provides a solution} of~Theorem~3.1.

\demo{Proof of Theorem~{3.1}}
 Since $\delta_j\to 0$, $T$ is an almost staircase transformation.
Moreover, $\# C_n\to\infty$ and $\# C_n\le n$ for all $n\in\Bbb N$.
Hence $T$ is mixing \cite{Ad}.
Moreover, it is mixing of all orders since $T$ is of rank one \cite{Ka}, \cite{Ry1}.
We also note that
$$
F_{s_1+\cdots+s_j}=F_0^{(j+1)}=F_{s_j}^{(j)}.
$$
Consider a mapping $\Phi_j$ which  sends a vector
$1_{[f]_{s_j}}\otimes v\in L^2(X_j,\mu_j)\otimes l^2(K_{s_j})$ to
a vector $1_{[f]_{s_1+\cdots+s_j}}\otimes v\in L^2(X, \mu) \otimes l^2(K_{s_j})$, $f\in F_{s_j}^{(j)}$, $v\in l^2(K_{s_j})$.
 Then $\Phi_j$ extends naturally to a linear isomorphism of $\Cal H_{s_j}^{(j)}$ onto the subspace
$$
\Cal H_{s_1+\cdots+s_j}:=\bigoplus_{f\in F_{s_1+\cdots+s_j}}1_{[f]_{s_1+\cdots+s_j}}\otimes l^2(K_{s_j})\subset L^2(X\times K,\mu\times\lambda_K)
$$
and $\|\Phi_jv\|=\sqrt{\frac{\mu(X_{s_1+\cdots+s_j})}{\mu_j(X_{s_j})}}\|v\|_j$ for all $v\in \Cal H_{s_j}^{(j)}$.
Moreover,
$$
\Phi U^q_{(T_j)_{\beta^{(j)}},\chi}v=U^q_{T_\beta,\chi}\Phi v
$$
whenever $v,U^q_{(T_j)_{\gamma^{(j)}},\chi}v\in\Cal H_{s_j}^{(j)}$.
Therefore for each $j>0$,
there exist   $Q_j>0$ and  a vector    $v_j\in\Cal H_{s_1+\cdots+s_j}\oplus\Cal H_{s_1+\cdots+s_j}$ such that
for each $w\in \Cal H_{s_1+\cdots+s_{j-1}}\oplus\Cal H_{s_1+\cdots+s_{j-1}}$,
$$
\inf_{a_{-q},\dots,a_q\in\Bbb R}\bigg\| w-\sum_{|q|\le Q_j}a_{q}
\bigg(U^q_{T_{\gamma},\chi}\oplus U^q_{T_{\gamma},\xi}\bigg)v_j\bigg\|\le \delta_j
\|w\|
$$
whenever  $\chi\not\sim\xi\in\widehat A$.
This yields \thetag{3-5}.
As we have already mentioned, \thetag{3-5} implies $\Cal M(T_{\gamma,H})=E$ and that $T_\gamma$ is weakly mixing.
Since $T$ is mixing of all orders, it follows that  $T_\gamma$ is mixing of all orders \cite{Ru}.
  Hence $T_{\gamma, H}$ is mixing of all orders because it is a factor of $T$. \qed
\enddemo

\head 4. Mixing transformations~whose~set~of~spectral multiplicities contains 2
\endhead

Let $E$ be an arbitrary  subset of \, $\Bbb N$
such that $2\in E$.
Our main purpose in this section is to prove the following theorem.

\proclaim{Theorem 4.1} There exists  a mixing (of all orders) transformation $S$ such that $\Cal M(S)=E$.
\endproclaim

We first establish  a generalization of the Algebraic lemma from \cite{KL} (cf. with \cite{Da5}).

\proclaim{Lemma 4.2}
Let $P$ be a subset of $\Bbb N$. Then there is a triplet $(K,B,D)$ such that $L(K,B,D)=P$. Moreover, if  $P$ is finite then we can choose $B$ finite and $K$  finite and cyclic.
\endproclaim
\demo{Proof}
Let $P=\{p_1,p_2,\dots,\}$ with $p_1<p_2<\cdots$.
For each $i>0$, we select a finite group $B_i$ and an automorphism $\theta_i$ of $B_i$ such that for each $0\ne b\in B_i$, the length of the $\theta_i$-orbit of $b$ is $p_i$.
We now set
$$
\align
B&:=B_1\oplus B_2^{\oplus p_1}\oplus B_3^{\oplus (p_1p_2)}\cdots,\\
D&:=\{(g_1,g_2,\dots)\in B\mid g_i=0\text{ if }i\notin\{1,2,2+p_1,2+p_1+p_1p_2,\dots\}\}.
\endalign
$$
We consider two automorphism of $B_i^{\oplus(p_1\cdots p_{i-1})}$.
The first one, say $\sigma_i$, is generated by the cyclic permutation of the coordinates.
The second one, say $\theta_i'$, is the Cartesian product $\theta_i\times\text{id}\times\cdots\times\text{id}$.
We note that $(\theta_i'\sigma_i)^{p_1\cdots p_{i-1}}=\theta_i\times\cdots\times\theta_i$.
We now define an automorphism $\theta$ of $B$ by setting
$$
\theta:=\theta_1\times \theta_2'\sigma_2\times \theta_3'\sigma_3\times\cdots.
$$
We consider $B$ equipped with $\theta$ as a $\Bbb Z$-module.
Let us compute the set $L(\Bbb Z,B,D)$.
Take $0\ne g=(g_1,g_2,\dots,)\in D$.
Let $l$ the the maximal non-zero coordinate of $g$.
If $l>2$ then
$l=2+\sum_{i=1}^kp_1\cdots p_k$ for some $k>0$.
It is easy to see that
$$
\min\{i>1\mid\text{such that }\theta^ig\in D\}=p_1\cdots p_{k+1}.
$$
Moreover, $\theta^{p_1\cdots p_{k+1}}g=(g_1,\dots,g_{l-1}, \theta_{k+2}g_l,0,0,\dots)$.
Hence the length of the $\theta$-orbit of $g$ is $p_{k+2}$.
The remaining cases, when $l=1,2$, are considered in a similar way.
We thus obtain $L(\Bbb Z,B,D)=P$.
Passing to a compactification as in the proof of Lemma~3.1, we obtain a compact Abelian group $K$ such that $B$ is a $K$-module and $L(K,B,D)=P$.
\qed
\enddemo

From now on and till the end of the section  $A$ denotes the dual group $\widehat B$ and $H$ denotes the dual group $\widehat{B/ D}\subset A$.
As follows from the proof of Lemma~4.2  we may assume without loss of generality that the conditions (i)--(iii) of Lemma~3.1 hold now as well.
We restate and prove Theorem~4.1 in the following way.

\proclaim{Theorem 4.3} There exists  a mixing (of all orders) transformation $S$ such that $\Cal M(S)=L(K,\widehat A,\widehat{A/H})\cup\{2\}$.
\endproclaim

We now state two auxiliary lemmata.

\proclaim{Lemma 4.4 \text{(\cite{Ag1}, \cite{Ry3})}}
Let $V$ be a unitary operator with a simple continuous spectrum.
\roster
\item"\rom{(i)}"
If the weak closure of powers of $V$ contains $\delta(I+V)$ for some $\delta>0$ then $V\odot V$ has a simple spectrum.
\item"\rom{(ii)}"
If $V\odot V$ has a simple spectrum then $V\otimes V$ has  a homogeneous spectrum of multiplicity~$2$.
\endroster
\endproclaim

\proclaim{Lemma 4.5 \text{(\cite{KaL})}} Let $V_i$, $i=1,2$, be unitary operators with a simple spectrum. Assume moreover that
\roster
\item"\rom{(i)}"
$V_i^{n_t}\to \delta(I+V_i^*)$ weakly, $i=1,2$, and
\item"\rom{(ii)}"
$V_i^{m_t}\to \delta(I+c_iV_i^*)$ weakly, $i=1,2$
\endroster
for some $\delta>0$.
If $c_1\ne c_2$ then $V_1\otimes V_2$ has also a simple spectrum.
\endproclaim

Let $0<i_n\le r_n$, $r_n\to\infty$,  $i_n/r_n\to\delta$, as $n\to\infty$, where $0<\delta< 1$.
Consider a partition
$$
\Bbb N=\bigsqcup_{a\in A_0}(\Cal N_a\sqcup\Cal M_a)\sqcup
\bigsqcup_{k\in K_0}\Cal L_k
$$
of $\Bbb N$ into infinite subsets $\Cal N_a$, $\Cal M_a$, $\Cal L_{k}$.

If $n+1\in \bigsqcup_{a\in A_0}\Cal N_a\sqcup
\bigsqcup_{k\in K}\Cal L_k$ then we define $C_{n+1},\alpha_{n+1},\beta_{n+1}$ by exactly the same formulae   as
in Section~3.
If $n+1\in\Cal M_a$ for some $a\in A_0$,
 we  let
$$
c_{n+1}(i):=
\cases
0, &\text{ if \ } i=0,\\
c_{n+1}(i-1)+h_n, &\text{ if \ } 0<i< i_n,\\
c_{n+1}(i-1)+h_n+1, &\text{ if \ } i_n\le i< 2i_n,\\
c_{n+1}(i-1)+h_n+i-2i_n, &\text{ if \ } 2i_n\le i< r_n,
\endcases
$$
 $C_{n+1}:=\{c_{n+1}(i)\mid 0\le i<r_n\}$,
$$
\alpha_{n+1}(c_{n+1}(i)):=
\cases
0, &\text{ if } i=0,\dots, i_n-1,\\
\alpha_{n+1}(c_{n+1}(i-1))- a, &\text{ if }i=i_n,\dots,2i_n-1,\\
\alpha_{n+1}(c_{n+1}(i-1)), &\text{ if }i=2i_n,\dots,r_n-1.
\endcases
$$
and  $\beta_{n+1}\equiv 1$ on $C_{n+1}$.

Let $(X,\mu,T)$ stand for  the $(C,F)$-system associated with $(C_{n+1},F_n)_{n\ge 0}$. Denote by $\gamma=(\beta,\alpha):\Cal R\to G$ the $(C,F)$-cocycle associated with $(\beta_n,\alpha_n)_{n>0}$.

\proclaim{Proposition 4.6} For each $a\in A_0$,
\roster
\item"\rom{(i)}" $U^{h_n}_{T_\beta,1}\to\delta \,I+(1-\delta)P_0$ weakly as $\Cal N_a-1\ni n\to\infty$, where $P_0$ is the orthogonal projection onto the subspace of constant functions on $X\times K$.
\item"\rom{(ii)}"
$U^{h_n}_{T_\beta,\chi}\to\delta \, l_\chi(a)I$ weakly as $\Cal N_a-1\ni n\to\infty$ for each $1\ne \chi\in\widehat A$.
\item"\rom{(iii)}" $U^{h_n}_{T_\beta,1}\to\delta(I+U^{*}_{T_\gamma,1} )+(1-2\delta)P_0$ weakly as $\Cal M_a-1\ni n\to\infty$, where $P_0$ is the orthogonal projection onto the subspace of constant functions on $X\times K$.
\item"\rom{(iv)}"
$U^{h_n}_{T_\beta,\chi}\to\delta( I+ l_\chi(a)U_{T_\gamma,\chi}^*)$ weakly as $\Cal M_a-1\ni n\to\infty$ for each $1\ne \chi\in\widehat A$.
\endroster
For each $k\in K$,
\roster
\item"\rom{(v)}" $U^{h_n}_{T,1}\to\delta \,I+(1-\delta)P_0$ weakly as $\Cal L_k-1\ni n\to\infty$, where $P_0$ is the orthogonal projection onto the subspace of constant functions on $X$.
\item"\rom{(vi)}"
$U^{h_n}_{T,\eta}\to\delta \, \eta(k)I$ weakly as $\Cal L_k-1\ni n\to\infty$ for each $1\ne \eta\in\widehat K$.
\endroster
\endproclaim

\demo{Proof}
Repeating the argument from the proof of Proposition~3.2  verbally
we obtain (i), (ii), (v), (vi).
The remaining  claims can be shown in  a similar way.
For instance, take $n\in\Cal M_a-1$,  subsets $B\subset F_n$ and $\eta\in\widehat K$. Then we obtain three ``sums'' (instead of the two sums in the proof of Proposition~3.2) in the righthandside
of the following:
 $$
\aligned
(U^{h_n}_{T,\eta} 1_{[B]_n})(x)
&=
\eta(\beta(x,T^{h_n}x))1_{T^{-h_n}[B]_n}(x)\\
&=\sum_{j=1}^{i_n-1}1_{[B+c_{n+1}(j-1)]_{n+1}}(x)\\
&+\sum_{j=i_n}^{2i_n-1}U^{*}_{T,\eta}1_{[B+c_{n+1}(j-1)]_{n+1}}(x)\\
&+
\sum_{j=2i_n+1}^{r_n-1}  U^{2i_n-j}_{T,\eta} 1_{[B+c_{n+1}(j-1)]_{n+1}}(x) +\bar o(1)
\endaligned
$$
for all $x\in X$ and $k\in K$.
Hence
$$
U^{h_n}_{T,\eta}\to \delta (I+U^{*}_{T,\eta})
+(1-2\delta)P_\eta \text{ \ weakly as }
\Cal M_a-1\ni n\to\infty.
$$
The claim~(iii) follows.
Similar minor ``changes'' of the proof of Proposition~3.2 are needed to establish (iv).
We leave details to the readers.
\qed
\enddemo

We now prove a {\it partially rigid} version of Theorem~4.3.

\proclaim{Theorem 4.7} $\Cal M(T_{\gamma,H}\times T_\beta)=L(K,\widehat A,\widehat{A/H})\cup\{2\}$.
\endproclaim

\demo{Proof}
Since $U_{T_{\gamma,H}}=\bigoplus_{\chi\in \widehat{A/H}}U_{T_\beta,\chi}$,
it is enough to show the following:
\roster
\item"(a)"
$U_{T_\beta,\chi}\otimes U_{T_\beta}$ has a simple spectrum for each $1\ne\chi\in\widehat A$,
\item"(b)" $U_{T_\beta}\otimes U_{T_\beta}$ has a homogeneous spectrum of multiplicity 2
    (in the orthocomplement to the constants),
\item"(c)"
$U_{T_\beta,\chi}\otimes U_{T_\beta}$ and $U_{T_\beta,\chi'}\otimes U_{T_\beta}$
are unitarily equivalent if $\chi'\sim\chi$,
\item"(d)"
the measures of maximal spectral type of $U_{T_\beta,\chi}\otimes U_{T_\beta}$ and $U_{T_\beta,\chi'}\otimes U_{T_\beta}$
are  mutually singular if $\chi'\not\sim\chi$.
\endroster

It follows from (v) and (vi) of Proposition~4.6 that $U_{T_{\beta},\chi}$
has a simple spectrum for each $\chi\in\widehat A$ (it was shown in the
proof of Theorem~3.7). Denote by $U_{T_\beta}^\circ$ the restriction of $U_{T_\beta}$ to
the orthocomplement to the constants. Then making use of (iii) and (iv) of
Proposition~4.6 we deduce from Lemma~4.5 that $U_{T_\beta,\chi}\otimes
U_{T_\beta}^\circ$ has a simple spectrum if $\chi\ne 1$. On the other hand, it
follows from (i) and (ii) of Proposition~4.6 that
$$
(U_{T_\beta,\chi}\otimes U_{T_\beta}^\circ)^{h_n}\to\delta^2 I
\text{ \ and \ }
U_{T_\beta,\chi}^{h_n}\to\delta I\text{ as }\Cal N_0-1\ni n\to\infty.
\tag4-1
$$
Therefore
the maximal spectral types of $U_{T_\beta,\chi}\otimes U_T^\circ$ and
$U_{T_\beta,\chi}$ are mutually disjoint. Hence the orthogonal sum of
these two operators, i.e. $U_{T_\beta,\chi}\otimes U_T$, has a simple
spectrum. Thus~(a) is proved.

We deduce from  Lemma~4.4 and Proposition~4.6(iii) that $U_{T_\beta}^\circ\otimes U_{T_\beta}^\circ$ has a homogeneous spectrum of multiplicity 2.
In view of Proposition~4.6(i),
$$
(U_{T_\beta}^\circ\otimes U_{T_\beta}^\circ)^{h_n}\to\delta^2 I\text{ \  and \ }
(U_{T_\beta}^\circ)^{h_n}\to\delta I\text{ as  }\Cal N_0-1\ni n\to\infty,
$$
 it follows that $U_{T_\beta}^\circ\otimes U_{T_\beta}^\circ$  and $U_{T_\beta}^\circ$ are spectrally disjoint. Therefore we obtain (b).

(c) follows from \thetag{1-1}.

To show (d), we first find $a\in A_0$ such that $l_\chi(a)\ne l_{\chi'}(a)$. It follows from (i) and (ii) of Proposition~4.6 that
$$
(U_{T_\beta,\chi}\otimes U_{T_\beta}^\circ)^{h_n}\to\delta^2 l_\chi(a) I
\text{ \ and \ }
(U_{T_\beta,\chi'}\otimes U_{T_\beta}^\circ)^{h_n}\to\delta^2 l_{\chi'}(a) I
\text{ as }\Cal N_a-1\ni n\to\infty.
$$
Hence $U_{T_\beta,\chi}\otimes U_{T_\beta}^\circ$ and $U_{T_\beta,\chi'}\otimes U_{T_\beta}^\circ$ are spectrally disjoint.
From \thetag{4-1} we deduce that $U_{T_\beta,\chi}\otimes U_{T_\beta}^\circ$ and $U_{T_\beta,\chi'}$ are spectrally disjoint as well as $U_{T_\beta,\chi'}\otimes U_{T_\beta}^\circ$ and $U_{T_\beta,\chi}$ are. Hence (d) follows.
 \qed
\enddemo

\demo{Proof of Theorem 4.3}
We note that the  condition (a) and (d) from the proof of Theorem~4.7 together are equivalent to the following
$$
(U_{T_\beta,\chi}\otimes U_{T_\beta})\oplus (U_{T_\beta,\xi}\otimes U_{T_\beta}) \text{ has a simple spectrum for all  $\chi\not\sim\xi \in\widehat A$}.\tag 4-2
$$
In view of Lemma~4.4(ii), the condition (b) of Theorem~4.7 is equivalent to
$$
U_{T_\beta}\oplus(U_{T_\beta}^\circ\odot U_{T_\beta}^\circ)\text{ has a simple spectrum.}\tag4-3
$$
  Thus \thetag{4-2} and \thetag{4-3} yield  that $\Cal M(T_{\gamma,H}\times T_\beta)=L(K,\widehat A,\widehat{A/H})\cup\{2\}$ and
$T_{\gamma}$  is weakly mixing.

It remains to  force mixing in the same way as we did in the proof of Theorem~3.1. For that we need to construct an {\it approximating} sequence $(T_j,\gamma^{(j)})_{j>0}$ of $\delta_j$-pairs with $\delta_j\to 0$ and pass to the $(C,F)$-limit. Our main concern now is to preserve \thetag{4-2} and \thetag{4-3} in the limit.
This requires only a slight modification of the argument in the proof of Theorem~3.3. We leave details to the reader.
\qed
\enddemo

\head 5. Concluding remarks
\endhead

\roster
\item
Modifying slightly the techniques developed in this paper one can obtain  other series of sets admitting mixing realizations.
For instance, consider the spectral multiplicities of the Cartesian products
$T_{\gamma,H}\times T_\beta^{\times r}$ with $r>1$ and their natural factors.
\item
Other mixing realizations may appear as Poissonian  (and Gaussian) ones.
For instance, given any $q>2$, the geometric sequence $\{q,q^2,q^3,\dots\}$ admits a mixing Poissonian realization \cite{DaR2}.
\item
We do not know, however, are there mixing realizations for $n\cdot E$, where $1\in E\subset \Bbb N$ and $n>2$.
Weakly mixing realizations for such subsets were constructed in \cite{Da3}.
 \endroster

\enddocument